\documentclass[onefignum,onetabnum]{siamonline171218}
\newcommand{\sech}{\ensuremath{\mathrm{sech}}}

\usepackage{lipsum}
\usepackage{amsfonts}
\usepackage{graphicx}
\usepackage{epstopdf}
\usepackage{algorithmic}
\ifpdf
  \DeclareGraphicsExtensions{.eps,.pdf,.png,.jpg}
\else
  \DeclareGraphicsExtensions{.eps}
\fi

\usepackage{enumitem}
\setlist[enumerate]{leftmargin=.5in}
\setlist[itemize]{leftmargin=.5in}

\newsiamremark{remark}{Remark}
\newsiamremark{hypothesis}{Hypothesis}
\crefname{hypothesis}{Hypothesis}{Hypotheses}
\newsiamthm{claim}{Claim}

\headers{Busse balloon deformation and splitting}{E. Siero and H. Uecker}

\title{Busse balloon deformation and splitting by non-local interaction: the influence of grazing on a Klausmeier vegetation model}

\author{E. Siero\thanks{Wageningen University 
  (\email{eric.siero@wur.nl}).}
\and H. Uecker\thanks{University of Oldenburg (\email{hannes.uecker@uni-oldenburg.de}).\\ June 5, 2026}}

\usepackage{amsopn}

\makeatletter
\newcommand*{\addFileDependency}[1]{
  \typeout{(#1)}
  \@addtofilelist{#1}
  \IfFileExists{#1}{}{\typeout{No file #1.}}
}
\makeatother

\ifpdf
\hypersetup{
  pdftitle={Busse balloon deformation},
  pdfauthor={E. Siero and H. Uecker}
}
\fi
\def\fnl{f_{\rm nl}}

\newcommand{\spr}[1]{\langle #1 \rangle}
\def\huga#1{\begin{gather} #1 \end{gather}}
\def\ddn{\frac{{\rm d}}{{\rm d}n}}
\def\pa{\partial}
\newcommand{\barr}{\begin{array}}\newcommand{\earr}{\end{array}}
\newcommand{\reff}[1]{(\ref{#1})}
\newcommand{\bpm}{\begin{pmatrix}}\newcommand{\epm}{\end{pmatrix}}

\def\al{\alpha}\def\R{{\mathbb R}}\def\ddt{\frac{{\rm d}}{{\rm d}t}}
\usepackage{hyperref} 
\usepackage{cleveref} 
\let\cref\ref
\let\Cref\ref
\usepackage{pstricks}

\def\BB{BB}
\def\BBs{BBs}
\newcommand{\bci}{\begin{itemize}}\newcommand{\eci}{\end{itemize}}

\begin{document}

\maketitle

\begin{abstract}
Large areas on all continents except Antarctica
are covered by dryland vegetation patterns, with wavelengths typically 
in the range of tens of meters. In models, many wavelengths are simultaneously
stable, and we argue that this multi-stability also holds for real world patterns. 
We then study the shape of the Busse balloon representation of multi-stability for a previously introduced Klausmeier model with non-local grazing. For
this we first extend application of singular perturbation to
Klausmeier/Gray-Scott models to include non-local interaction,
providing analytical control near the homoclinic limit, and then 
use numerical continuation to demonstrate deformation of Busse
balloons away from the typical banana shape, in four non-local grazing
regimes. Since the Busse balloon has been invoked to support ``evasion
of tipping'', we underscore the importance of the shape of the Busse
balloon when inferring ecosystem response to, e.g., climate change.
\end{abstract}

\begin{keywords}
  pattern formation, numerical continuation, Busse balloon, dryland vegetation, Turing patterns
\end{keywords}

\begin{AMS}
  35B36, 35B25, 35B60, 92D40
\end{AMS}

\section{Introduction}
In pattern forming systems, patterns with different wavelengths can 
simultaneously be stable. This multi-stability is conveniently represented in a so-called Busse balloon (BB), 
pioneered by F. Busse in the study 
of a Rayleigh-B\'enard experiment \cite{Busse1978}, where a fluid is heated from below. As the temperature forcing is increased, heat transport by conduction is replaced by regular convection patterns that subsequently disappear in the turbulent regime. In repeated physical experimentation, 
the number of convection rolls that appear is not strictly reproducible, 
but rather is from a range of admissible values \cite{BusseWhitehead1971}. 
This was corroborated by simultaneous stability 
(i.e.\,for identical parameter values) of the patterns in the mathematical 
fluid mechanics model. 

In drylands, vegetation patterns are ubiquitous \cite{Deblauweetal2008}. After publication of the Klausmeier model in 1999 \cite{Klausmeier1999}, reaction-(advection-)diffusion systems became popular models for their study (e.g. also \cite{HilleRisLambersetal2001,Giladetal2004}). In these models, as the annual precipitation drops below a critical Turing instability \cite{Turing1952} value, homogeneous vegetation is replaced by spatial patterns through self-organization. These patterns subsequently disappear if precipitation drops below a value where any vegetation can be sustained.

The concept of the \BB\ was introduced in spatial ecology in the past 
two decades \cite{VdSteltetal2013,Sherratt2013a,Siteuretal2014}. 
Here multi-stability of different pattern configurations suggests 
that past differences may lead to different resulting patterns at 
equal environmental conditions (parameter values) \cite{Sherratt2015}. Analysis of remote sensing imagery at two locations in Somalia have shown a large spread in dryland vegetation pattern wavenumbers which could not be explained by variations in precipitation or terrain slope, in agreement with the existence of a \BB\ in the real ecosystem \cite{Bastiaansenetal2018}. However, the possibility could not be fully excluded that variables outside scope (e.g.\,differences in soil composition) explain wavenumber spread, so that selection of a unique wavenumber (mono-stability) is a priori still possible. 

Vegetation patterns are viewed as a means of the ecosystem to mitigate drought stress, allowing vegetation survival at precipitation levels that would otherwise be too low (reviewed in \cite{Rietkerketal2021}). Another driver of desertification is grazing. In many models, grazing is incorporated in a local mortality term, ignoring the mobility of grazers \cite{Siero2018}. In this paper we show how the shape of the \BB\ changes in four different grazing regimes. 

The region with low wavenumber and low precipitation admits a
treatment using singular perturbation
(e.g.\,\cite{DoelmanEckhausKaper2000}), exploiting the scale
separation of the diffusion coefficients of water and vegetation.
Near this homoclinic limit the profile of the vegetation
component has the form of a pulse, which for the Klausmeier
model has an explicit approximate expression. We find algebraic
relations for the locus of the fold bifurcation, and thus the
existence boundary, for the Klausmeier model with non-local grazing
terms. Near the homoclinic limit, the existence boundary and stability
boundary lie asymptotically close (e.g. \cite{VdSteltetal2013}), so
this analysis can be used to infer the shape of the \BB\ in
this corner at low computational cost.

The paper is organized as follows. In the remainder of
the introduction we present the Klausmeier model, including four types
of non-local grazing. In Section \cref{sec:generic} we explain why observation
of a trend in wavenumbers is clear evidence for
multi-stability, so that \BBs\ are generically
applicable. In Section \cref{sec:pulsefold} we determine equations for
approximation of the fold curve at low precipitation $a$ (where
pulse-type patterns disappear), extending application of singular
perturbation to Klausmeier/Gray-Scott type systems with non-local
terms. In Section \cref{sec:changingshape} we present the \BBs\ in
the four different grazing regimes based on continuation of the full
PDE system, including a comparison with the approximation of the fold
curves based on Section \cref{sec:pulsefold}. An interesting feature is that 
in two regimes \BBs\ may become disconnected, yielding stability 
of patterns at high and low precipitations, and/or at low and high 
wavenumbers, but not in the repective intermediate ranges. 
A discussion follows in Section \cref{sec:discussion}.

\subsection{Dryland vegetation Klausmeier model}
We use the extended Klausmeier model \cite{Klausmeier1999,Siteuretal2014} as a phenomenological model for spatial vegetation pattern formation by self-organization. It is a spatially explicit two-component 
(surface water $w$, vegetation $n$) model, which we 
treat on a large one-dimensional interval 
with periodic boundary conditions. 
In scaled form, on flat terrain, the Klausmeier model is given by 
\begin{equation}
	\left\{\begin{alignedat}{5} \label{eq:extKG}
		\frac{\pa w}{\pa t}=&&d\frac{\pa^2w}{\pa x^2}&&+a&&-w&&-wn^2, \\
		\frac{\pa n}{\pa t}=&&\frac{\pa^2n}{\pa x^2}&&-g_jn^j&&-m_0n&&+wn^2. 
 \end{alignedat}\right.
\end{equation} The terms $d\frac{\pa^2w}{\pa x^2}$ and $\frac{\pa^2n}{\pa x^2}$ model water diffusion and plant dispersal, notably $d\gg 1$. The parameter $a$ models precipitation; $-w$ models evaporation; the terms $\pm wn^2$ model water uptake by the vegetation and $m_0$ is vegetation mortality unrelated to grazing. The additional term $g_jn^j$ representing grazing mortality is explained below (Section \cref{sec:grazing}). 

The non-linear dependence of water uptake $wn^2$ on vegetation $n$ is a crucial ingredient of the Klausmeier model, and often referred to as the cause of local facilitation. Formally, such a non-linear term can be derived by reduction of a three-component model (with both surface and soil water), where both water infiltration and uptake depend linearly on vegetation, to a two-component (surface water only) model \cite{Siero2020}. We fix parameters in accordance with those used in \cite{Siteuretal2014}: water diffusion coefficient $d=500$ and non-grazing death rate $m_0=0.45$; the precipitation parameter $a$ is varied.

\subsubsection{Including grazing} \label{sec:grazing}
To understand how grazing alters the desertification process, we add grazing terms to the extended Klausmeier model. Grazing involves a varying distribution of herbivores, with variable feeding demands depending on a variable amount of forage. Based on fundamental theory on herbivore movement and preferences (from \cite{Holling1959a,FretwellLucas1969,Noy-Meir1975}, among others), a non-local and 
non-linear grazing mortality term can be derived. Here we only present a minimal explanation of the resulting grazing terms, and refer the reader to \cite{Siero2018,Sieroetal2019} for modeling details.

A measure of forage in the domain $[-L,L]$ is given by $\langle n^j\rangle=\frac{1}{2L}\int_{-L}^Ln^j(x)dx$. For the aggregation parameter $j$ we make a distinction between two cases. The case $j=1$ is referred to as \emph{proportional} grazing; here the number of herbivores at a location is proportional to amount of biomass at this location. In the other case $j=2$, referred to as 
\emph{disproportionate} grazing, higher stacks of biomass attract 
disproportionately many herbivores.
The dependence of grazing pressure on overall forage also comes in two flavors. The case where the grazing pressure is a monotonically decreasing function of forage is referred to as \emph{sustained} grazing. The case where it is initially an increasing function of forage and then switches to decreasing, is referred to as \emph{natural} grazing. Combining proportional/disproportionate with sustained/natural grazing, we arrive at the following four grazing pressure functions:
\begin{equation} \label{eq:grazing}
	\left\{\begin{alignedat}{4}
		&g_{1,\mathrm{sus}}(\langle n\rangle)=\frac{m}{K+\langle n\rangle} &&\textnormal{(proportional sustained);} && \\
		&g_{1,\mathrm{nat}}(\langle n\rangle)=\frac{m\langle n\rangle}{{K}^2+{\langle n\rangle}^2} &&\textnormal{(proportional natural);} \\
		&g_{2,\mathrm{sus}}(\langle n^2\rangle)
=\frac{m}{K+\langle n^2\rangle}&&\textnormal{(disproportionate sustained);} \\
		&g_{2,\mathrm{nat}}(\langle n^2\rangle)=\frac{m\langle n^2\rangle}{{K}^2+{\langle n^2\rangle}^2}\qquad&&\textnormal{(disproportionate natural).} \\
 \end{alignedat}\right.
\end{equation}

\begin{figure}
 \centering
 \begin{minipage}{0.49\textwidth}
	\centering sustained grazing \\
  \includegraphics[width=0.99\textwidth, clip=true, trim=130 335 145 260]{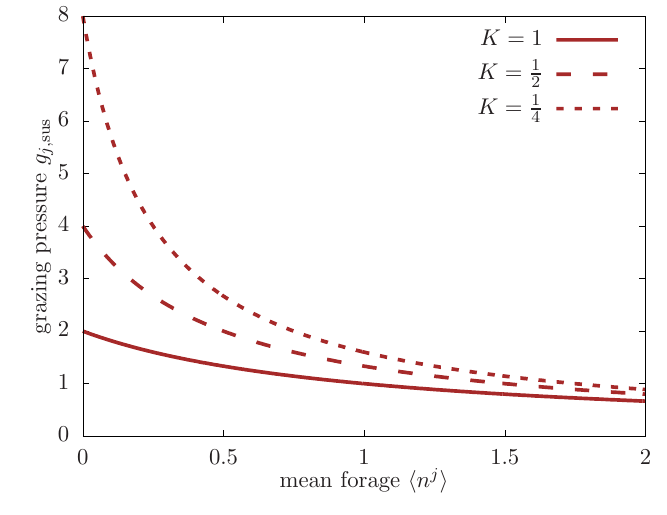}
 \end{minipage}
 \begin{minipage}{0.49\textwidth}
	\centering natural grazing \\
  \includegraphics[width=0.99\textwidth, clip=true, trim=130 335 145 260]{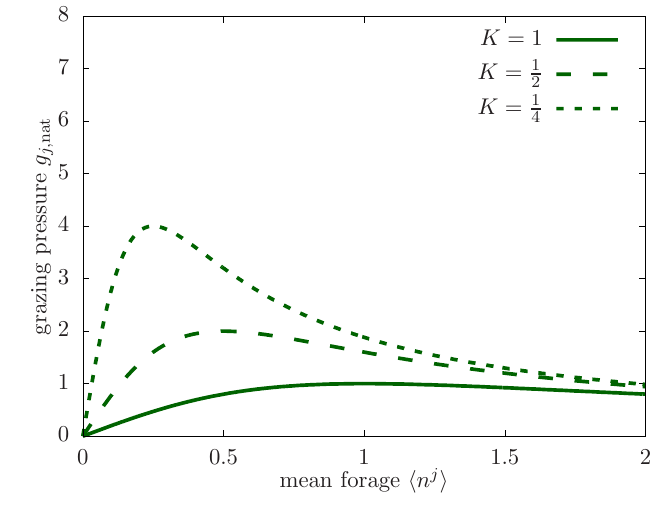}
 \end{minipage}
	\caption{Sustained and natural grazing pressure as a function of mean forage $\langle n^j\rangle$, for $m=2$ and $K=1,\frac{1}{2}$,$\frac{1}{4}$. In both cases the grazing pressure $g_j$ is approximately inversely proportional to mean forage in case of abundant forage. The function $g_{j,\mathrm{sus}}$ starts at $\frac{2}{K}$ and then monotonically decreases; in contrast $g_{j,\mathrm{nat}}$ starts at $0$ and first initially increases to a maximum of $\frac{1}{K}$ at $\langle n^j\rangle=K$.} \label{fig:modeling}
\end{figure}

The parameter $m$ is the maximal mean grazing rate, the spatial average of the grazing rate in case of superfluous forage availability, 
i.e.\,$m=\lim_{\langle n^j\rangle\rightarrow\infty} \langle g_jn^j\rangle$. 
The parameter $K$, referred to as half-persistence, is the amount of forage $\langle n^j\rangle$ where the mean grazing rate is half the maximal mean grazing rate, i.e.\,$\langle g_j(K)n^j\rangle=\frac{m}{2}$. 
For the grazing pressure,
$g_j(K)=\frac{m}{2K}$ is half the maximum value in case of sustained grazing and the maximum for natural grazing (see Fig. \cref{fig:modeling}). The lower the value of $K$, the more persistent the grazing. Here we vary between $K=1,\frac{1}{2}$ or $\frac{1}{4}$. Multiple values for $m$ have also been used, but 
we fix $m=2$, since variation of $K$ yielded the most interesting results. 

\section{Genericity of multi-stability of patterns and the Busse balloon} \label{sec:generic} 
The perceived multi-stability of vegetation patterns has been a source
of debate. On the one hand, observation of a wide spread in wave
numbers for a given set of parameters \cite{Bastiaansenetal2018}
supports the idea. On the other hand, mono-stability with noisy
variation could arguably also result in a spread. Here we present a
simple model independent argument explaining why mono-stability of
real world patterns is not robust against parameter changes. For
illustration, and as a baseline for the subsequent including of grazing, 
we first briefly turn to the Klausmeier model without grazing and its
\BB.

\begin{figure}[H]
 \includegraphics{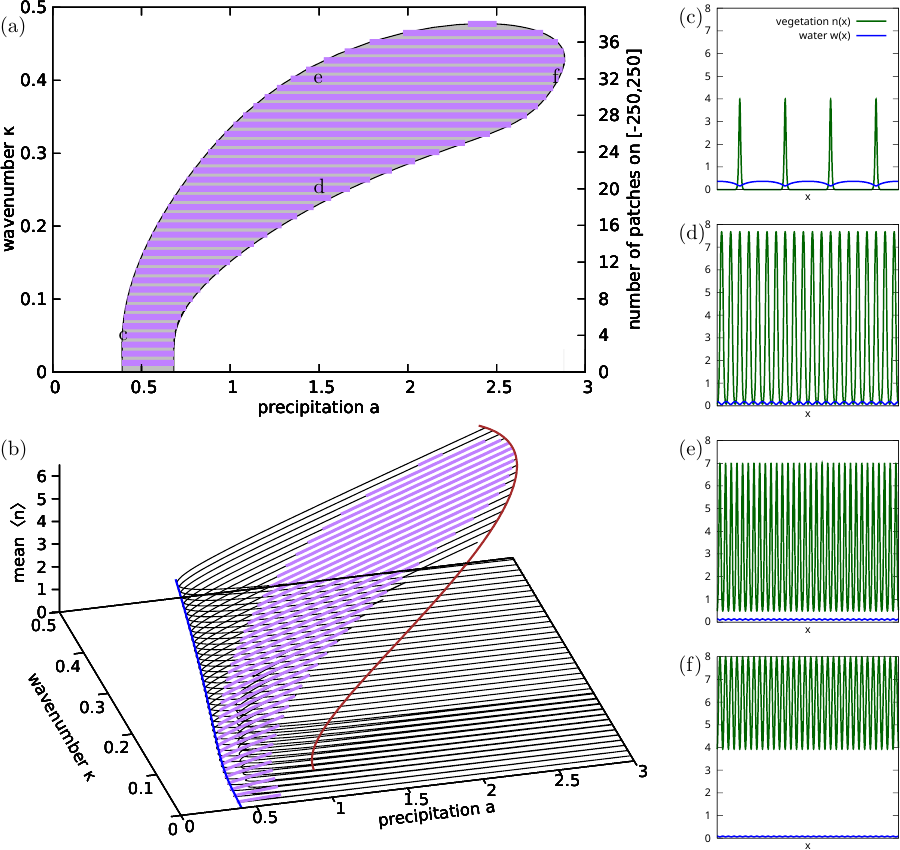}
	\caption{\BB\ for the Klausmeier model without grazing. 
(a) Black curve: sideband instability curve (from \cite{Siteuretal2014}) 
delimiting the \BB\ on the unbounded domain. 
Horizontal purple lines: representation of the \BB\ 
on the bounded domain $[-250,250]$, 
corresponding to integer numbers of patches. 
Letters c-f correspond to patterns in panels (c)-(f). 
(b) Purple lines from (a) with an additional axis of 
mean vegetation $\langle n\rangle$, and black curved lines representing 
unstable patterns; brown curve is the Turing instability line, 
and the blue curve (see \eqref{eq:foldnograzing}) approximates a fold of 
pulsed patterns. (c) Pattern profile with four 
patches ($\kappa\approx 0.05$) at $a=0.4$; 
(d) $20$ patches ($\kappa\approx 0.25$) at $a=1.5$; 
(e) $32$ patches ($\kappa\approx 0.4$) at $a=1.5$; 
(f) $32$ patches ($\kappa\approx 0.4$) at $a=2.8$. 
(d),(e) Two patterns with distinct wavenumber that are stable at the same parameter values. 
(e),(f) Two patterns with the same wavenumber at distinct precipitation $a$. 
}  \label{fig:Bb0}
\end{figure}

The \BB\ is a representation of multi-stability of spatially
periodic patterns. Fig. \cref{fig:Bb0}(a) shows the \BB\ 
representation of stable spatially periodic vegetation patterns of the
Klausmeier model without grazing, both for the bounded domain
$[-250,250]$ and the unbounded domain $(-\infty,\infty)$. On the
unbounded domain, for a single value of the precipitation parameter
$a$, there is a continuous range of wavenumbers $\kappa$ inside the
balloon, each representing a stable periodic vegetation pattern. On
the bounded domain, patterns are required to satisfy the periodic
boundary conditions, so the range of wavenumbers is quantized to a
discrete set of values.

Now suppose that a trend in real world pattern wavenumber is observed, for instance that wavenumber increases with precipitation (e.g. \cite{Deblauweetal2011}). Together with mono-stability this would imply that the wavenumber is a non-constant function of precipitation. We will first discuss the case of 
bounded domains (which is the most relevant in the real world) but also extend the argument to the unbounded domain.

\paragraph{Bounded domains}
On a bounded domain, only a discrete set of wavenumbers is possible, since the boundary conditions need to be satisfied. Mono-stability of patterns on an interval $[a_1,a_2]$ of precipitation values would mean that for every $a\in [a_1,a_2]$ there exists a unique periodic pattern with a certain wavenumber $\kappa$. Because of the observation of a trend, constant wavenumber is not possible, so there must be a value of $a=a_\mathrm{edge}$ that corresponds to the stability boundary of two patterns with distinct wavenumbers. These patterns are not asymptotically close to each other. Now, by changing any parameter value other than precipitation (e.g. mortality), the stability boundary $a_\mathrm{edge}$ of both wavenumbers will shift. Generically, this shift will be unequal for the two wavenumbers, creating either a gap or overlap in stability. Thus, mono-stability does not persist.

\paragraph{Unbounded domain}
On the unbounded domain, there are no boundary conditions that make the set of admissible wavenumbers discrete. Importantly, wavenumbers act distinctly 
from e.g.\,amplitudes or spatial averages, in the sense that a small 
change in wavenumber never results in a pattern that is close to the original pattern, because somewhere on the unbounded domain the patterns will be in anti-phase. Put differently, a small change in wavenumber can only occur by passing through the set of aperiodic functions (disregarding growing domains). So, even on the unbounded domain with a continuous range of wavenumbers available, the wavenumber cannot change in a continuous way, and we can use the same arguments as in the bounded domain case to explain that mono-stability is not robust.

In summary, observation of a trend in pattern wavenumber is sufficient empirical evidence to establish multi-stability, which in models is 
characterized by the \BB.

\section{Fold approximation using singular perturbation} \label{sec:pulsefold}
Pulse-like patterns like in Fig. \cref{fig:Bb0}(c) exist at low
precipitation values $a$, and we can use their specific structure to
approximate a fold curve by continuation of three algebraic equations, 
which can then be compared to the
fold captured by continuation of the full PDE model. For
this, we will now construct a single pulse solution (as in
e.g. \cite{DoelmanEckhausKaper2000}) at the center of the domain
$[-L,L]$, so that the wavenumber of the corresponding pulse pattern is
$\kappa=\frac{L}{\pi}$.
Water diffusion is much faster than vegetation dispersal ($d\gg 1$),
so that $\epsilon:=\frac{1}{\sqrt{d}}$ with $d=500$ makes
$\epsilon\approx 0.045$ a small parameter. Using this, the Klausmeier
model \eqref{eq:extKG} can be viewed as a singularly perturbed
system. We will incorporate both the proportional and disproportionate
grazing, represented by $g_1n$ and $g_2n^2$ respectively, in one
go. 

Multiplying the $w$-equation of \eqref{eq:extKG} by
$\epsilon^2$, we obtain
\begin{equation} \label{eq:extKG2}
 \left\{\begin{aligned} 
  \epsilon^2w_t=&w_{xx}+\epsilon^2(a-w-wn^2)\ , \\
	 n_t=&n_{xx}-(m_0+g_1)n+(w-g_2)n^2\ .
 \end{aligned}\right. 
\end{equation} 
After setting $\epsilon=0$, the inner solution (near the origin) is given by 
\begin{equation} \label{eq:sech}
 \left\{\begin{alignedat}{1}
	 w(x)=&w_0\ , \\
	 n(x)=&\frac{3(m_0+g_1)}{2(w_0-g_2)}\sech^2\left(\frac{\sqrt{m_0+g_1}x}{2}\right)\ .
 \end{alignedat}\right. 
\end{equation} On the other hand, after rescaling $\xi=\epsilon x$, \eqref{eq:extKG} becomes
\begin{equation} \label{eq:extKG3}
 \left\{\begin{aligned} 
  w_t=&w_{\xi\xi}+a-w-wn^2, \\
	 n_t=&\epsilon^2n_{\xi\xi}-(m_0+g_1)n+(w-g_2)n^2, 
 \end{aligned}\right. 
\end{equation} which after setting $\epsilon=0$ yields the outer solution 
\begin{equation}
 \left\{\begin{alignedat}{1}
	n(\xi)=&0, \\
	 w(\xi)=&a-(a-w_0)\frac{\cosh(|\xi|-\epsilon L)}{\cosh\left(\epsilon L\right)}\ .
 \end{alignedat}\right. 
\end{equation} 
By Fenichel theory \cite{Fenichel1979}, for $\epsilon$ sufficiently small there exists a steady state near the combined pulse-type solution
\begin{equation} \label{eq:pulse}
 \left\{\begin{aligned}
	 w_p(x)=&a-(a-w_0)\frac{\cosh(|\epsilon x|-\epsilon L)}{\cosh\left(\epsilon L\right)}, \\
 n_p(x)=&\frac{3(m_0+g_1)}{2(w_0-g_2)}\sech^2\left(\frac{\sqrt{m_0+g_1}x}{2}\right).
 \end{aligned}\right. 
\end{equation} 
However, $w_0$ with $g_2<w_0<a$ still needs to de determined, given the 
two conservation laws 
\begin{align} 
		\int_{-L}^L a-w_p-w_pn_p^2\ dx=&0, \label{eq:conservation1} \\
		\int_{-L}^L -(m_0+g_1)n_p+(w_p-g_2)n_p^2\ dx=&0. \label{eq:conservation2}
\end{align} As a preparation, we compute 
\begin{align}
	\int_{-L}^Lw_pdx=&
	2aL-\frac{2(a-w_0)\tanh\left(\epsilon L\right)}{\epsilon}; \label{eq:avwp} \\
	\int_{-L}^Ln_p(x)dx=&\frac{3(m_0+g_1)}{2(w_0-g_2)}\int_{-L}^L\sech^2\left(\frac{\sqrt{m_0+g_1}}{2}x\right)dx \label{eq:avnp} \\
    =&\frac{6\sqrt{m_0+g_1}}{w_0-g_2}\tanh\left(\frac{\sqrt{m_0+g_1}L}{2}\right); \notag \\
	\int_{-L}^Ln_p^2(x)dx=&\frac{9(m_0+g_1)^2}{4(w_0-g_2)^2}\int_{-L}^L\sech^4\left(\frac{\sqrt{m_0+g_1}}{2}x\right)dx \label{eq:avnp2} \\ 
	=&\frac{3(m_0+g_1)^\frac{3}{2}}{(w_0-g_2)^2}\left(\sech^2\left(\frac{\sqrt{m_0+g_1}L}{2}\right)+2\right)\tanh\left(\frac{\sqrt{m_0+g_1}L}{2}\right). \notag
\end{align} 
For $L\gg 1$, from \eqref{eq:avnp} and \eqref{eq:avnp2} we obtain: 
\begin{align}
	\int_{-L}^Ln_p(x)dx\approx \frac{6\sqrt{m_0+g_1}}{w_0-g_2},&&\textnormal{ so }&&\langle n_p\rangle:=&\frac{1}{2L}\int_{-L}^Ln_p(x)dx\approx \frac{3\sqrt{m_0+g_1}}{L(w_0-g_2)}; \label{eq:foragej1} \\
	\int_{-L}^Ln_p^2(x)dx\approx \frac{6(m_0+g_1)^\frac{3}{2}}{(w_0-g_2)^2},&&\textnormal{ so }&&\langle n_p^2\rangle:=&\frac{1}{2L}\int_{-L}^Ln_p^2(x)dx\approx \frac{3(m_0+g_1)^\frac{3}{2}}{L(w_0-g_2)^2}. \label{eq:foragej2}
\end{align}

Since $w_p\approx w_0$ for $n_p\not\approx 0$, using \eqref{eq:avwp} and \eqref{eq:foragej2}, the left-hand side of \eqref{eq:conservation1} can be approximated by
\begin{equation}
		\int_{-L}^L a-w_p-w_pn_p^2\ dx\approx 2aL-\left(2aL-\frac{2(a-w_0)\tanh(\epsilon L)}{\epsilon}\right)-w_0\frac{6(m_0+g_1)^\frac{3}{2}}{(w_0-g_2)^2}, \\
\end{equation} 
so the first conservation law \eqref{eq:conservation1} boils down to 
\begin{equation} \label{eq:cont}
	f(a,w_0):=(a-w_0)(w_0-g_2)^2\tanh(\epsilon L)-3\epsilon w_0(m_0+g_1)^\frac{3}{2}=0.
\end{equation}

Again since $w_p\approx w_0$ for $n_p\not\approx 0$, now using \eqref{eq:foragej1} and \eqref{eq:foragej2}, the left-hand side of \eqref{eq:conservation2} can be approximated by
\begin{equation}
		\int_{-L}^L -(m_0+g_1)n_p+(w_p-g_2)n_p^2\ dx\approx -\frac{6(m_0+g_1)^\frac{3}{2}}{w_0-g_2}+\frac{6(m_0+g_1)^\frac{3}{2}}{w_0-g_2}=0,
\end{equation} 
so \eqref{eq:conservation2} is always (approximately) fulfilled.

For no grazing, i.e., $g_1=g_2=0$, \eqref{eq:cont} simplifies to
\begin{equation} \label{eq:conservationwaternograzing}
	f(a,w_0)=(a-w_0)w_0\tanh(\epsilon L)-3\epsilon m_0^\frac{3}{2}=0.
\end{equation} We identify the fold of pulses as the locus where the implicit function theorem breaks down, i.e., 
\begin{equation}
	\frac{\pa f}{\pa w_0}=(a-2w_0)\tanh(\epsilon L)=0,
\end{equation} which implies $w_0=\frac{a}{2}$, so that with \eqref{eq:conservationwaternograzing} we obtain
\begin{equation} \label{eq:foldnograzing}
	a=2\sqrt{\frac{3\epsilon m_0^\frac{3}{2}}{\tanh(\epsilon L)}}=2\sqrt{\frac{3\epsilon m_0^\frac{3}{2}}{\tanh\left(\frac{\epsilon\pi}{\kappa}\right)}},
\end{equation} 
which is used to compute the blue fold curve in Fig. \cref{fig:Bb0}(b).

If we do incorporate grazing, it is not possible to find $a$ as an explicit function of $\kappa$. To compute a solution branch with the unknowns $a$, $\kappa$ (or equivalently $L$), $w_0$, and $\langle n^j\rangle$, we establish three algebraic equations linking these unknowns. 
We first distinguish between the case of proportional ($j=1$) and disproportionate ($j=2$) grazing and for each of these subsequently discriminate between sustained and proportional grazing. 

\subsection{Fold continuation for proportional grazing}\label{sec:foldprop}
For proportional grazing, $g_2=0$, so \eqref{eq:cont} simplifies to
\begin{equation} \label{eq:conservationwaterj1}
	f(a,w_0)=(a-w_0)w_0\tanh(\epsilon L)-3\epsilon (m_0+g_1)^\frac{3}{2}=0,
\end{equation} and the implicit function theorem breaks down if
\begin{equation} \label{eq:Dw0conservationwaterj1}
	\frac{\pa f}{\pa w_0}=(a-2w_0)\tanh(\epsilon L)-\frac{9}{2}\epsilon\sqrt{m_0+g_1}\frac{\pa g_1}{\pa w_0}=0,
\end{equation} and from \eqref{eq:foragej1} we obtain 
\begin{align} 
	\frac{\pa \langle n_p\rangle}{\pa w_0}\approx &\frac{3}{L}\left(\frac{-1}{w_0^2}\sqrt{m_0+g_1}+\frac{1}{2w_0}\frac{1}{\sqrt{m_0+g_1}}\frac{\pa g_1}{\pa w_0}\right)  \label{eq:Dw0foragej1}\\
    =&\frac{-3\kappa}{2\pi w_0^2}\frac{1}{\sqrt{m_0+g_1}}\left(2(m_0+g_1)-w_0\frac{\pa g_1}{\pa w_0}\right).\notag 
\end{align}

For sustained proportional grazing, $g_1=\frac{m}{K+\langle n\rangle}$, so
\begin{equation}
	\frac{\pa g_1}{\pa w_0}=\frac{-m}{\left(K+\langle n\rangle\right)^2}\frac{\pa \langle n\rangle}{\pa w_0}, 
\end{equation} 
which combined with \eqref{eq:Dw0foragej1} 
yields
\begin{equation} \label{eq:Dw0gj1sus}
	\frac{\pa g_1}{\pa w_0}\approx\frac{6\kappa m(m_0+g_1)}{2\pi w_0^2\sqrt{m_0+g_1}(K+\langle n\rangle)^2+3\kappa mw_0}.
\end{equation} 
\eqref{eq:foragej1}, \eqref{eq:conservationwaterj1} and \eqref{eq:Dw0conservationwaterj1} with \eqref{eq:Dw0gj1sus} are used to compute the blue fold curve in Fig. \cref{fig:Bbj1sus}(b).

For natural proportional grazing $g_1=\frac{m\langle n\rangle}{K^2+\langle n\rangle^2}$, so
\begin{equation}
	\frac{\pa g_1}{\pa w_0}=\frac{m(K^2-\langle n\rangle^2)}{(K^2+\langle n\rangle^2)^2}\frac{\pa \langle n\rangle}{\pa w_0}, 
\end{equation} 
and combined with \eqref{eq:Dw0foragej1} 
this yields
\begin{equation} \label{eq:Dw0gj1nat}
	\frac{\pa g_1}{\pa w_0}\approx\frac{6\kappa m(m_0+g_1)(\langle n\rangle^2-K^2)}{2\pi w_0^2\sqrt{m_0+g_1}(K^2+\langle n\rangle^2)^2+3\kappa mw_0(\langle n\rangle^2-K^2)}. 
\end{equation}
\eqref{eq:foragej1}, \eqref{eq:conservationwaterj1} and \eqref{eq:Dw0conservationwaterj1} with \eqref{eq:Dw0gj1nat} are used to compute the blue fold curve in Fig. \cref{fig:Bbj1nat}(b).

\subsection{Fold continuation for disproportionate grazing}
\label{sec:folddis}
For $g_1=0$, \eqref{eq:cont} simplifies to
\begin{equation} \label{eq:conservationwaterj2}
	f(a,w_0)=(a-w_0)(w_0-g_2)^2\tanh(\epsilon L)-3\epsilon w_0m_0^\frac{3}{2}=0,
\end{equation} 
and the implicit function theorem breaks down if
\begin{align} 
		\frac{\pa f}{\pa w_0}=&\left(-(w_0-g_2)^2+2(a-w_0)(w_0-g_2)\left(1-\frac{\pa g_2}{\pa w_0}\right)\right)\tanh(\epsilon L)-3\epsilon m_0^\frac{3}{2} \label{eq:Dw0conservationwaterj2} \\
	=&(w_0-g_2)\left(-(w_0-g_2)+2(a-w_0)\left(1-\frac{\pa g_2}{\pa w_0}\right)\right)\tanh(\epsilon L)-3\epsilon m_0^\frac{3}{2}=0, \notag
\end{align} 
and from \eqref{eq:foragej2} we obtain 
\begin{equation} \label{eq:Dw0foragej2}
	\frac{\pa \langle n_p^2\rangle}{\pa w_0}\approx\frac{3}{L}\cdot\frac{-2m_0^\frac{3}{2}}{(w_0-g_2)^3}\left(1-\frac{\pa g_2}{\pa w_0}\right)=\frac{6\kappa}{\pi}\frac{m_0^\frac{3}{2}}{(w_0-g_2)^3}\left(\frac{\pa g_2}{\pa w_0}-1\right).
\end{equation}

For sustained disproportionate grazing, $g_2=\frac{m}{K+\langle n^2\rangle}$, hence 
\begin{equation}
	\frac{\pa g_2}{\pa w_0}=\frac{-m}{\left(K+\langle n^2\rangle\right)^2}\frac{\pa \langle n^2\rangle}{\pa w_0}, 
\end{equation} 
which combined with \eqref{eq:Dw0foragej2} 
yields
\begin{equation} \label{eq:Dw0gj2sus}
	\frac{\pa g_2}{\pa w_0}\approx\frac{6\kappa mm_0^\frac{3}{2}}{\pi (w_0-g_2)^3(K+\langle n^2\rangle)^2+6\kappa mm_0^\frac{3}{2}},
\end{equation} 
and \eqref{eq:foragej2}, \eqref{eq:conservationwaterj2} and \eqref{eq:Dw0conservationwaterj2} with \eqref{eq:Dw0gj2sus} are used to compute the blue fold curve in Fig. \cref{fig:Bbj2sus}(b).

For natural disproportionate grazing, 
$g_2=\frac{m\langle n^2\rangle}{K^2+\langle n^2\rangle^2}$, hence 
\begin{equation}
	\frac{\pa g_2}{\pa w_0}=\frac{m(K^2-\langle n^2\rangle^2)}{(K^2+\langle n^2\rangle^2)^2}\frac{\pa \langle n^2\rangle}{\pa w_0}, 
\end{equation} 
which combined with \eqref{eq:Dw0foragej2} 
yields
\begin{equation} \label{eq:Dw0gj2nat}
	\frac{\pa g_2}{\pa w_0}\approx\frac{6\kappa mm_0^\frac{3}{2}(\langle n^2\rangle^2-K^2)}{\pi (w_0-g_2)^3(K^2+\langle n^2\rangle^2)^2+6\kappa mm_0^\frac{3}{2}(\langle n^2\rangle^2-K^2)}, 
\end{equation} 
and \eqref{eq:foragej2}, \eqref{eq:conservationwaterj2} and \eqref{eq:Dw0conservationwaterj2} with \eqref{eq:Dw0gj2nat} are used to compute the blue fold curve in Fig. \cref{fig:Bbj2nat}(b).

\section{Shape changes and splittings of the Busse balloons}\label{sec:changingshape}
We now numerically investigate the shape of the \BB\ in each of the four 
grazing regimes \eqref{eq:grazing}, in each case for the 
half-persistence parameter $K=1$, $\frac{1}{2}$ and $\frac{1}{4}$. 
We combine: 
\bci
\item[1)] computation of the Turing instability where the patterns emerge (simple relationship between $a$ and $\kappa$, Appendix \cref{app:Tur}); 
\item[2)] approximation of the fold where pulse patterns disappear for small precipitation $a$ (continuation of three algebraic equations (Section \cref{sec:pulsefold}), using \emph{auto} \cite{auto}); 
\item[3)] numerical continuation of the full PDE model including spectral stability (using \emph{pde2path}, Section \cref{sec:Bbcont}).
\eci  
Thus, we effectively take two--dimensional slices through the three-dimensional 
(parameters $a,K$, and wavenumber $\kappa$) 
\BBs, and naturally, we may expect these 
slices to change continously in $K$. An interesting feature is that 
for disproportinate grazing these slices become disconnected at small 
$K$, leaving no stable patterns at some intermediate $a$ and $\kappa$. 
Additionally, we shall focus on the shapes of the \BBs\ at small $a$ and 
$\kappa$ and the consequences for tipping and possible desertification, 
and on the shapes of the \BBs\ at large $a$ and small $\kappa$, yielding 
some rather unexpected multistability between homogeneous vegetation 
and large wavelength patterns.  



\subsection{Proportional sustained grazing}
\begin{figure}[ht]
 \includegraphics{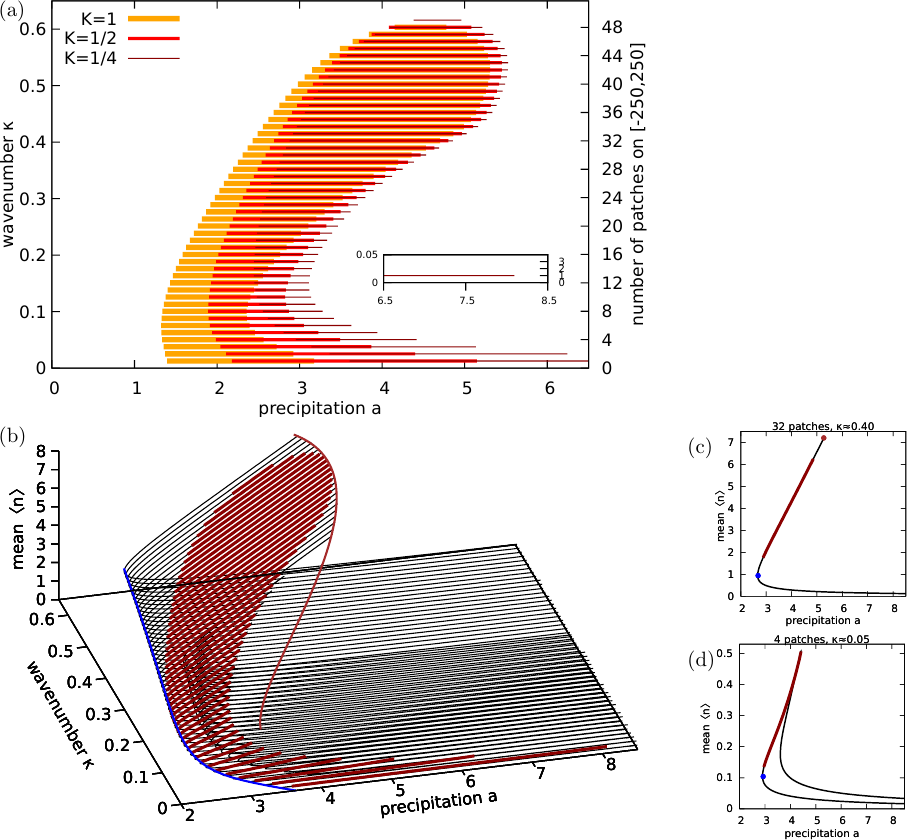}
	\caption{\BBs\ 
for proportional sustained grazing. (a) Horizontal orange, red and dark-red bars represent stable patterns for half-persistence $K=1$, $K=1/2$ and $K=1/4$ respectively. 
(b) \BB\ of $K=1/4$ from (a) with additional $\langle n\rangle$ axis,
and unstable patterns (black), the Turing curve (brown), 
and the approximate fold curve (Section \cref{sec:foldprop}, blue). 
(c) Slice with $32$ patches ($\kappa\approx 0.4$) 
attached to the Turing instability curve, which behaves as a pitchfork bifurcation: the solution flips phase and the continuation doubles back. (d) Slice with 
$4$ patches ($\kappa\approx 0.05$, detached from the Turing instability curve.}  \label{fig:Bbj1sus}
\end{figure}

The resulting \BBs\ for proportional sustained grazing are shown in Fig. \cref{fig:Bbj1sus}. In panel (a) we observe that in the bottom left corner, the \BB\ lifts off from the $\kappa{=}0$ axis, unlike the benchmark \BB\ without grazing (Fig. \cref{fig:Bb0}). This effect becomes stronger as the grazing half-persistence $K$ decreases, so when grazing becomes more persistent. Related to this, the stability boundary at the bottom left does not emerge vertically from the $\kappa=0$ axis, but tilts to the left. In this corner, small wavenumber combined with low precipitation implies a low mean forage, which for sustained grazing results in an unsustainable large grazing pressure (overgrazing), providing a heuristic explanation why (stable) patterns do not appear here. The smaller the half-persistence $K$ the more pronounced this effect is (confer Fig. \cref{fig:modeling}). Since the \BB\ constrains where the system can reside, at low precipitation $a$ this opens up the scenario of a transition from a relatively `healthy' large wavenumber pattern immediately to the bare desert state ($\kappa=0$), revealing limitations on the ability of pattern formation to mitigate drought in this scenario.

At the same time, in the bottom right corner the \BB\ is considerably extending to the right. For $K=1/4$, this results in a range of simultaneous stability of low wavenumber patterns and the homogeneously vegetated state (which is stable to the right of the brown Turing instability curve of Fig. \cref{fig:Bbj1sus}(b)). The singular perturbation approximation of the fold of pulse-type solutions (Section \cref{sec:foldprop}) in blue performs very well, even for larger wavenumbers.

\subsection{Proportional natural grazing}
\begin{figure}[ht]
  \includegraphics{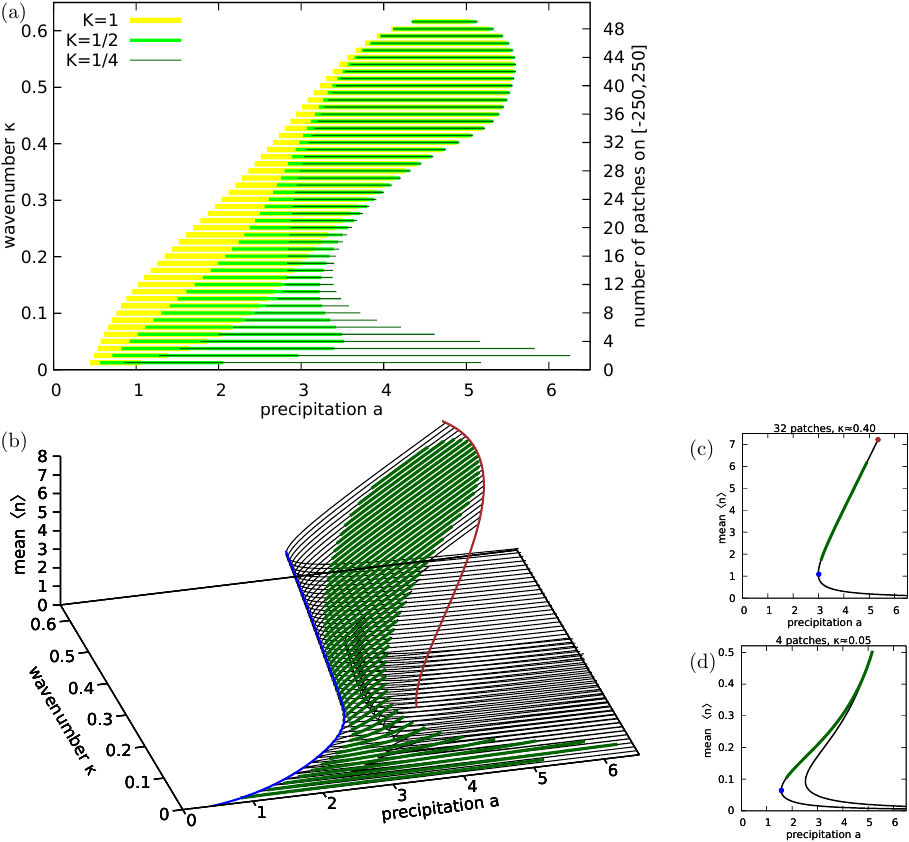}
	\caption{\BB\ 
for proportional natural 
grazing. (a) Horizontal yellow, green and dark-green bars represent stable patterns for half-persistence $K=1$, $K=1/2$ and $K=1/4$ respectively. 
(b) \BB\ of $K=1/4$ from (a) with additional $\langle n\rangle$ axis, and unstable patterns 
(black), the Turing curve (brown), 
and the approximate fold curve (Section \cref{sec:foldprop}, blue).  
(c) Slice with $32$ patches ($\kappa\approx 0.4$)  attached to the Turing instability curve. (d) Slice with 
$4$ patches ($\kappa\approx 0.05$, detached from the Turing instability curve.}  \label{fig:Bbj1nat}
\end{figure}

\BBs\ for proportional natural grazing are shown in Fig. \cref{fig:Bbj1nat}. Contrary to proportional sustained grazing, the bottom left corner in panel (a) now extends to the left and the left boundary is tilted to the right; these effects again become more pronounced as the grazing half-persistence $K$ decreases. This relates to the fact that now the grazing pressure is relatively low at low mean forage (Fig. \cref{fig:modeling}). Proportional natural grazing introduces a qualitative change halfway the \BB: the curvature of the upper boundary goes from convex to concave. As this progresses for smaller $K$, less 
patterns with intermediate wavenumbers are available. Natural grazing pressure is largest at the intermediate amount $\langle n\rangle=1/K$ of mean forage, providing a heuristic explanation for the pinching of the \BB\ halfway. Unavailability of intermediate wavenumbers creates conditions for a possible jump from high to low wavenumber at intermediate precipitation $a$. 

As in the case of proportional sustained grazing, the bottom right corner of the \BB\ extends to the right. For $K=1/4$ this again results in a range of simultaneous stability of the homogeneously vegetated state and low wavenumber patterns - but not for the pattern with a single patch in the domain. The approximation of the fold of pulse-type solutions by the blue curve based on singular perturbation (Section \cref{sec:foldprop}) is again surprisingly good. 

\subsection{Disproportionate sustained grazing}
\begin{figure}[ht]
 \includegraphics{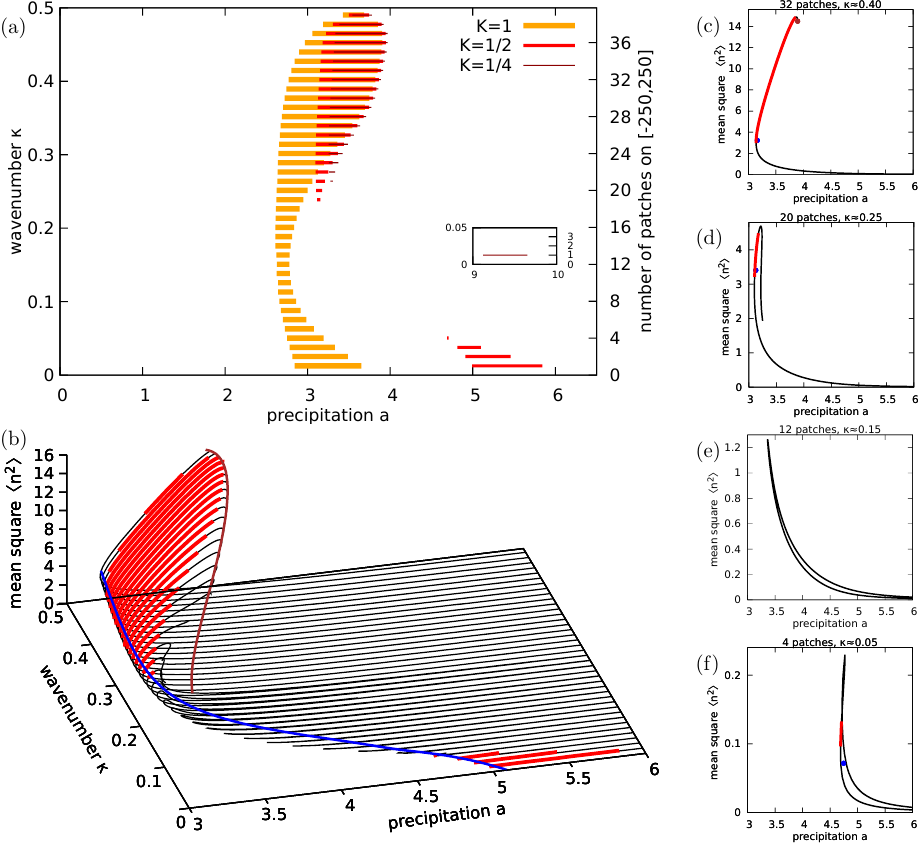}
	\caption{\BBs\ 
for disproportionate sustained grazing. (a) Horizontal orange, red and dark-red bars represent stable patterns for half-persistence $K=1$, $K=1/2$ and $K=1/4$ respectively.  
(b) \BB\ of $K=1/2$ from (a) with additional $\langle n\rangle$ axis,
and unstable patterns (black), the Turing curve (brown), 
and the approximate fold curve 
(Section \cref{sec:folddis}, blue). (c) Slice with 
$32$ patches ($\kappa\approx 0.4$), attached to the Turing instability curve. 
(d) $20$ patches ($\kappa\approx 0.25$), with a very narrow stability range. 
(e) $12$ patches ($\kappa\approx 0.15$), with empty stability range. 
(f) Slice with $4$ patches ($\kappa\approx 0.05$), which 
has regained a pair of folds and a stability range.}  \label{fig:Bbj2sus}
\end{figure}

The \BBs\ for disproportionate sustained grazing are shown in Fig. \cref{fig:Bbj2sus}(a). For the relatively high value $K=1$, the shape of the \BB\ is similar to that of proportional sustained grazing. As $K$ is decreased towards $K=1/2$ and $K=1/4$, the balloon splits up into two parts. 
On a heuristic level, the disproportionate aggregation of herbivores at peaks of vegetation negates pattern formation when the grazing pressure is high enough, 
i.e.\ for small half-persistence $K$ and when precipitation $a$ and wavenumber $\kappa$ are low. 

Small amplitude Turing patterns emerge in the top part of the \BB, as
shown by the connection to the brown curve in panel (b). Since for the
top part there are no alternative stable patterns with low wavenumber
at the same range of precipitation $a$, direct transitions from high
wavenumber to the bare desert state are conceivable. The bottom part
of the \BB\ consists of small wavenumber patterns at high
precipitation.

In this case the singular perturbation approximation of the fold of pulse-type solutions (Section \cref{sec:folddis}) in blue is not as good as in the other cases, especially between $\kappa=0.1$ and $\kappa=0.2$. The bifurcation diagram of the pattern with $12$ patches in panel (e) exemplifies this, with the blue dot corresponding to the blue curve missing because it is too far off. Here the pulse type solutions exhibit a double hump (pattern profile not shown) that conflicts with the single hump solution assumed in \eqref{eq:sech}. This coincides with the presence of only a single fold in panel (e), compared to multiple folds in panel (d) and (f).

\subsection{Disproportionate natural grazing}
\begin{figure}[ht]
 \includegraphics{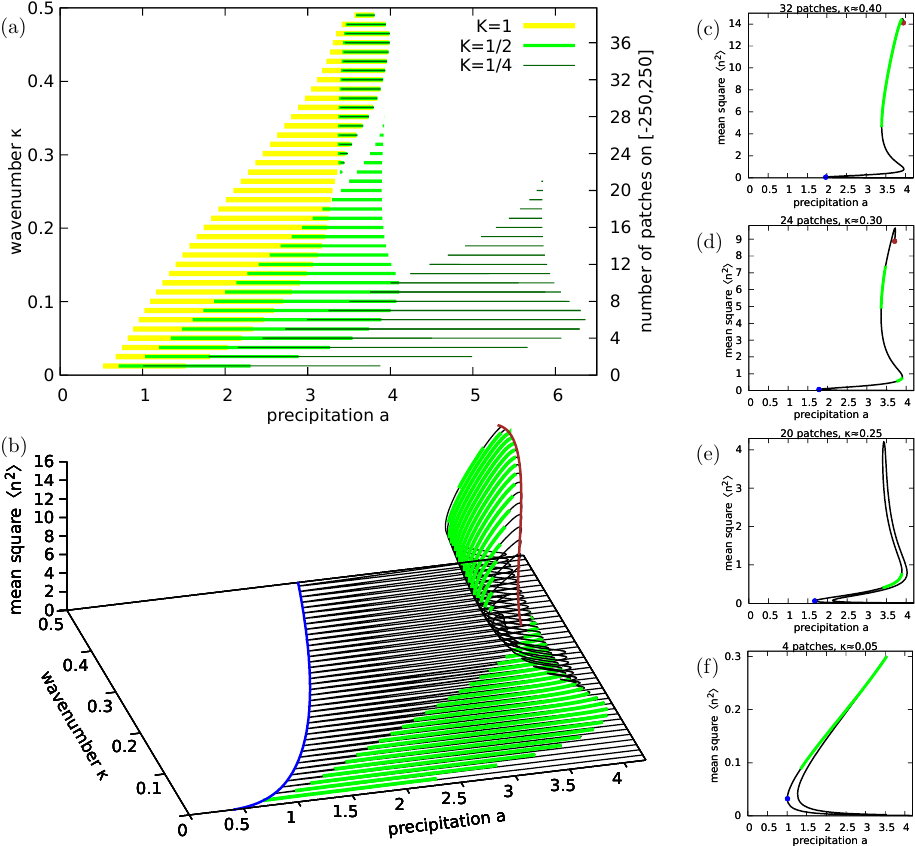}
	\caption{\BBs\ 
for disproportionate natural grazing (a) Horizontal yellow, green and dark-green bars represent stable patterns for half-persistence $K=1$, $K=1/2$ and $K=1/4$ respectively. \BBs\ for $K=1/2$ and $K=1/4$ have split. 
(b) \BB\ of $K=1/2$ from (a) with additional $\langle n\rangle$ axis,
and unstable patterns (black), the Turing curve (brown), 
and the approximate fold curve (Section \cref{sec:folddis}, blue). 
(c,d) Slices with $32$ patches ($\kappa\approx 0.4$) and  
$24$ patches ($\kappa\approx 0.3$), both attached to the 
Turing instability curve. In (d), 
an additional range of stability has appeared at low $\langle n^2\rangle$. 
(e,f) $20$ patches ($\kappa\approx 0.25$) and $4$ patches ($\kappa\approx 0.05$), without stable patterns at high $\langle n^2\rangle$. 
} \label{fig:Bbj2nat}
\end{figure}

For natural disproportionate grazing, the \BBs\ are presented in Fig. \cref{fig:Bbj2nat}(a). For $K=1$, the \BB\ is similar in shape to that of proportional natural grazing. As in the case of disproportionate sustained grazing, for $K$ decreasing towards $K=1/2$ and $K=1/4$ the balloon splits up into two parts, but here the bottom part extends from low to high precipitation $a$.

Small amplitude Turing patterns again emerge in the top part of the \BB, from the brown curve in panel (b). Since for the full range of precipitation values $a$ of the top part there are also stable patters in the bottom part, transition from high wavenumber to low wavenumber offers a viable alternative to 
transition from high wavenumber to the bare desert state directly. The singular perturbation approximation of the fold of pulse-type solutions (Section \cref{sec:folddis}) in blue this time is again spot on, even for large wavenumbers.

\section{Discussion} \label{sec:discussion}
Busse balloons provide an overview of stable patterns for a 
range of parameter values, and have made their way into ecology, 
also for plant communities \cite{Ferreetal2025} and beyond the study of dryland vegetation patterning \cite{BennettSherratt2019,Detmeretal2025}. 
In Section \cref{sec:generic} we provided arguments why observation of a trend in 
wavenumber - e.g.\,as a function of precipitation - is sufficient 
evidence to invalidate mono-stability, 
underscoring the importance of \BBs. These balloons are at the basis of 
the concept of tipping evasion \cite{Rietkerketal2021}, claiming that critical transitions can become gradual if spatial pattern formation is enabled. Such an evasion crucially depends on the size and shape of the \BB.

The theory of tipping evasion is inspired by typical banana shaped \BBs, 
as in Fig. \cref{fig:Bb0} \cite{Rietkerketal2021,Banerjeeetal2026}. This balloon has the property that the last patterns to destabilize before the desertification process is completed, have large wavelengths (see also \cite{Doelmanetal2018}, where this is referred to as Ni's conjecture). These localized patterns are near the homoclinic limit and thus permit analysis through singular perturbation (exploiting the fact that water diffuses faster than vegetation). Based on algebraic equations for the Klausmeier model with non-local grazing terms derived in Section \cref{sec:pulsefold}, we found that the way the fold curve is tilted depends on the type of grazing. In both cases with sustained grazing, the fold is tilted in such a way that the large wavelength patterns are {\em not} the last to become non-existent. Since in the homoclinic limit generally the existence and stability boundary are close by, this suggests that the large wavelength patterns are also not the last to destabilize. Computation of the entire \BB\ via continuation of the full system of PDEs confirms this (Fig. \cref{fig:Bbj1sus} and \cref{fig:Bbj2sus}).

It is hard to predict which wavenumbers in the \BB\ are selected 
\cite{Aschetal2025,EigentlerSensi2026}. Nevertheless, the \BB\ constrains complex dynamics and informs us which stable patterns are (un)available. As such, 
the tilting of the \BB\ near the homoclinic limit in the cases of 
sustained grazing makes a direct transition from intermediate wavelength 
to the bare desert state feasible. Thus, here pattern formation in response to water stress does not negate the possibility of a critical transition due to overgrazing. Similarly, the cases of natural grazing make a critical transition from high to low wavenumber possible, since intermediate wavenumbers are not always available (Fig. \cref{fig:Bbj1nat} and \cref{fig:Bbj2nat}).

Validating a mathematical model for a landscape scale ecological system is in many ways harder than that for a fluid mechanics lab setup \cite{Gandhietal2019}. Currently, time series with sufficient spatial resolution are still very short, but opportunities for machine learning mechanistic models from remote sensing data are opening up \cite{ZeeMarcosSieroPP}. Alternatively, it may be possible to optimize the matching of patterns observed over a wider area (the so-called empirical \BB, \cite{Bastiaansenetal2018}) with \BBs\ computed for a library of models, as a kind of space-for-time substitution. Since the \BB\ changes predictably when a parameter is varied (e.g. changes of 
Fig. \cref{fig:Bbj1sus}--\cref{fig:Bbj2nat} in response to $K$), an optimal match could be found iteratively. An extra challenge here is that not all stable patterns are equally likely to be observed, also in the light of stochasticity \cite{HamsterHeijsterSiero2025}.

\appendix
\section{Pattern formation: The Turing bifurcation curves} \label{app:Tur}
To compute where small amplitude patterns bifurcate, we first need to
compute the relevant spatially homogeneous steady state
$(\bar w,\bar n)$. From the first equation of \cref{eq:extKG} it
follows that $\bar w=\frac{a}{1+\bar n^2}$, so we only need to find
$\bar n$. The relevant $(\bar w,\bar n)$ turns out to always be the
one with maximal $\bar n>0$. From the second equation of
\cref{eq:extKG}, this means numerically solving
\begin{equation} \label{eq:barn} \left\{\begin{alignedat}{5}
  &\textnormal{prop.sust.:}&& 
m_0\bar n^3+(m_0K+m-a)\bar n^2+(m_0-aK)\bar n+m+m_0K&&=0, \\
 &\textnormal{prop.nat.:}&& 
m_0\bar n^4+(m-a)\bar n^3+(m_0+m_0K^2)\bar n^2+(m-aK^2)\bar n+m_0K^2&&=0, \\
 &\textnormal{disprop.sust.:}\quad && 
m_0\bar n^4+(m-a)\bar n^3+(m_0+m_0K)\bar n^2+(m-aK)\bar n+m_0K&&=0, \\
 &\textnormal{disprop.nat.:}&& 
m_0\bar n^6+(m-a)\bar n^5+m_0\bar n^4+m\bar n^3+m_0K^2\bar n^2-aK^2\bar n+m_0K^2&&=0,
 \end{alignedat}\right. 
\end{equation}
and selecting the largest solution, as a function of $a$. 

Next we need to linearize about $(\bar w,\bar n)$.
Since we want to determine Turing instability, we focus
on instability against spatially inhomogeneous perturbations, so we
restrict to perturbations $(\tilde n,\tilde w)$ with zero mean, so
$\langle \tilde n\rangle=\frac{1}{2L}\int_{-L}^L\tilde
n(x)\mathrm{d}x=0$ (and similarly for $\tilde w$). 
Following \cite{Siero2016}, the Gateaux derivative of $\langle n^j\rangle$ in
the $\tilde n$ direction is given by 
\begin{equation}
\begin{aligned}
 d(\langle n^j\rangle;\tilde n)=&
\lim_{h\rightarrow 0}\frac{\langle (n+h\tilde n)^j\rangle-\langle n^j\rangle}{h} 
 =\lim_{h\rightarrow 0}\frac{\int_{-L}^L\int_0^1\frac{d}{ds}(n+hs\tilde n)^jdsdx}{2Lh} \\
 =&\lim_{h\rightarrow 0}\frac{1}{2L}\int_{-L}^L\int_0^1\frac{d}{ds}j(n+hs\tilde n)^{j-1}
\tilde ndsdx =jn^{j-1}\langle\tilde n\rangle=0
\end{aligned}  
\end{equation} where we have ignored technical details for interchanging limit and integral. By the chain rule, for $g_{j,\mathrm{sus}}(\langle n^j\rangle)$ and $g_{j,\mathrm{nat}}(\langle n^j\rangle)$ (see \cref{eq:grazing}) it follows that 
\begin{equation}
d(g_{j,\mathrm{sus}}(\langle n^j\rangle);\tilde n)=
d(g_{j,\mathrm{nat}}(\langle n^j\rangle);\tilde n)=0. 
\end{equation} 
Thus, for perturbations with zero mean, the
linearization of proportional sustained grazing
$g_{1,\mathrm{sus}}(\langle n\rangle)n$ is
$\mathrm{graz}_\mathrm{lin}{=}g_{1,\mathrm{sus}}(\langle n\rangle)$,
that of proportional natural grazing
$g_{1,\mathrm{nat}}(\langle n\rangle)n$ is
$\mathrm{graz}_\mathrm{lin}=g_{1,\mathrm{nat}}(\langle n\rangle)$,
that of disproportionate sustained grazing
$g_{2,\mathrm{sus}}(\langle n\rangle)n^2$ is
$\mathrm{graz}_\mathrm{lin}=2g_{2,\mathrm{sus}}(\langle n\rangle)n$, 
and that of disproportionate natural grazing
$g_{2,\mathrm{nat}}(\langle n\rangle)n^2$ is
$\mathrm{graz}_\mathrm{lin}=2g_{2,\mathrm{nat}}(\langle n\rangle)n$. 
Thus, after linearizing all other terms in
\cref{eq:extKG} and replacing $\frac{\partial^2}{\partial x^2}$ by
$-k^2$, the dispersion relation that needs to be solved is 
\begin{equation} \label{eq:disp} 
\det\left(\begin{array}{cc}-dk^2-1-\bar n^2-\lambda & -2\bar w\bar n \\ 
\bar n^2 & -k^2-m_0-\mathrm{graz}_\mathrm{lin}+2\bar  w\bar n-\lambda
\end{array}\right)=0, 
\end{equation} with $\lambda=0$ at instability. 
The relationship between $a$ ($\bar n$ is a function of $a$ through \cref{eq:barn}) and 
$\kappa$ resulting from \cref{eq:disp} can easily be determined numerically.

\section{Busse balloons by numerical continuation} \label{sec:Bbcont}
Busse balloons were computed with the numerical continuation software
{\em pde2path} \cite{UeckerWetzelRademacher2014,p2pbook,p2phome}.   
For the domain we use the
interval $[-250,250]$ with $n_p=4096$ equidistant mesh points, and 
periodic boundary conditions and an associated standard translational 
phase condition.  The non-standard part here relates to the
nonlocal grazing term, which yields rank--1 corrections to otherwise 
sparse Jacobians, see also \cite[\S8.3.4]{p2pbook}. 
Here, writing the FEM formulation of \reff{eq:extKG} as 
$M\ddt \bpm w\\n\epm=G(w,n):=(G_1(w,n),G_2(w,n))^T$, where $M$ is the 
FEM mass matrix, we need to compute 
\huga{\label{Jform} 
J=\bpm \pa_w G_1&\pa_n G_1\\
\pa_w G_2&\pa_n G_2\epm=A-cb^T, 
}
where in the FEM discretization, $A\in\R^{2n_p\times 2n_p}$ comes from the 
Laplacian and the local terms and is sparse, 
while $cb^T$ with $c,b\in\R^{2n_p}$ is 
a full matrix but of rank 1. 
Writing the non--local terms $g_jn^j$ in $G_2$ 
as $\fnl(n,a), a=\spr{h(n)}$, i.e., 
\huga{
\fnl(n,a)=\left\{\barr{ll}  \frac m{K+a} n&\text{sustained,}\\
\frac {ma}{K^2+a^2} n&\text{natural,}\earr\right.
\text{ and }
h(n)=\left\{\barr{ll}n&\text{proportional,}\\ n^2&
\text{disproportionate,}\earr\right.
}
we have, similar to \reff{eq:disp}, 
\huga{\label{gcgj}
\left(\ddn\fnl(n,a)\right)\tilde n=\pa_n\fnl(n,a)\tilde n
+\pa_a\fnl(n,a)\spr{h'(n) \tilde n}. 
}
Thus, the $n$--parts $c_2,b_2$ of 
the vectors $c=(0,c_2)\in\R^{2n_p}$ and $b=(0,b_2)\in\R^{2n_p}$ are 
given by 
\huga{
c_2=\left\{\barr{ll}  -\frac m{(K+a)^2} n&\text{sustained,}\\
\frac {m(K^2-a^2)}{(K^2+a^2)^2} n&\text{natural,}\earr\right.
\text{ and }
b_2=\left\{\barr{ll}\al&\text{proportional,}\\ 
2\al\otimes n&\text{disproportionate,}\earr\right. 
}
where $\al\in\R^{n_p}$ is the weight vector for evaluating 
$\spr{n}$ as $\sum_{l=1}^{n_p} \al_l n_l$ with $n_l$ the nodal 
value of $n$ at point $l$, and $\al\otimes n$ denoting the 
componentwise product. 

Now, linear systems $Jz{=}G$ with $J$ from \reff{Jform} and $G{\in}\R^{2n_p}$, 
which appear in Newton loops, 
and also in the inverse vector iteration for computing selected small 
eigenvalues of $J$, can be treated efficiently by the 
Sherman--Morrison--Woodbury (SMW) formula \cite[\S2.7.3]{numrc} 
\huga{\label{smf} 
z=A^{-1}G+\al_0 
A^{-1}cb^TA^{-1}G, \quad \al_0=\frac 1{1-b^TA^{-1}c}.  
}
In this, and in particular in the inverse vector iteration for eigenvalues, 
we reuse approximate $LU$--decompositions of $A$ as much as possible.

For computation of the Busse balloon, the workflow now reads as follows:
\begin{enumerate}
\item Start with stable homogeneously vegetated state, decrease
precipitation and Turing bifurcate to patterned
state. 
\item Find lower and upper bound of stability
interval. 
\item Find solution at middle of stability interval. %
\item Jump to neighboring (admissible on finite domain) wavenumber.
\item Converge to stable pattern (this might fail) and return to step
2.
\end{enumerate}
Files to initialize generation of \BBs\ 
will appear at 
\cite[Tab ``Applications'']{p2phome}, together with some 
further technical description.



\section*{Acknowledgments}
ES was supported by a Humboldt postdoctoral fellowship in the initial phase of this study.

\bibliographystyle{plain}
\bibliography{bbbib}

\end{document}